\documentclass[12pt]{article}
\usepackage[francais]{babel}
\usepackage{amsmath,amsfonts,amssymb,amsthm}
\begin{document}

\title{Hasard et d\'eterminisme chez Laplace}
\author{Jean--Pierre kahane}
\date{}
\maketitle

\hfill\begin{minipage}{8cm}
Il n'y a rien  de contingent dans la nature des choses ; elles sont au contraire d\'etermin\'ees par la n\'ecessit\'e de la nature divine \`a exister et op\'erer d'une mani\`ere certaine.

Spinoza, Ethique I, proposition 29
\end{minipage}

\vskip1cm

Il est de bon ton aujourd'hui de traiter de haut le d\'eterminisme laplacien, comme une vision obtuse de la complexit\'e des choses. Certes le monde et la science ont beaucoup chang\'e en deux si\`ecles, mais la condamnation syst\'ematique de Laplace me parait un aspect de l'obscurantisme contemporain. Le d\'eterminisme de Laplace est fonction de son \'epoque, de la tradition cart\'esienne (voir la citation de Spinoza en exergue) et des progr\`es immenses que la m\'ecanique de Newton a fait r\'ealiser dans la compr\'ehension du monde. Il traduit les succ\`es de l'analyse math\'ematique comme m\'ethode  g\'en\'erale et universelle (c'est ainsi qu'il faut lire la phrase fameuse, reproduite ci-dessous, sur ``une intelligence qui, pour un instant donn\'e, conna\^{i}trait toutes les forces dont la nature est anim\'ee et la situation respective des \^etres qui la composent...'' ; il ne lui suffit pas de tout conna\^{i}tre, il lui faut encore \^etre ``assez vaste pour soumettre ces donn\'ees \`a l'analyse'', entendez, l'analyse math\'ematique). Il se rattache \`a une vision corpusculaire de la mati\`ere, que nous devons remettre en question avec la m\'ecanique quantique. Mais il me para\^{i}t en avance sur son temps, et m\^eme sur le n\^otre, par sa mani\`ere d'articuler d\'eterminisme, certitude et probabilit\'es. Avec Laplace, la notion de hasard, au lieu de se r\'eduire \`a l'inexplicable contingence, devient un auxiliaire de la science. Ce qu'il dit de la place souhaitable de la th\'eorie des probabilit\'es, \`a ``faire entrer dans le syst\`eme de l'instruction publique'', n'a vu ses premi\`eres r\'ealisations en France qu'il y a 40 ans, et reste d'actualit\'e.

Cet article est b\^ati autour de citations de Laplace, assez \'etendues pour bien rendre compte de sa pens\'ee. Il se ressent \'evidemment  de l'admiration sans r\'eserve que j'\'eprouve pour Laplace comme math\'ematicien. Le texte de r\'ef\'erence est l'``Essai philosophique sur les probabilit\'es'', de 1814.

\section*{Laplace math\'ematicien}

Laplace est un tr\`es grand nom des math\'ematiques. Ses \oe uvres, en quatorze gros volumes, se trouvent dans certaines biblioth\`eques. Il vaut la peine d'y aller voir : on y trouve des tr\'esors. C'est \`a la fois un virtuose du calcul, un cr\'eateur de concepts, et un ma\^{i}tre de la langue. Je m'attarderai ici sur deux \oe uvres \'ecrites \`a l'intention des non--sp\'ecialistes : les dix le\c cons profess\'ees aux Ecoles Normales en 1795, et l'essai philosophique sur les probabilit\'es (1814), d\'eveloppement de la dixi\`eme le\c con de 1795, qui eut un tel succ\`es, du vivant de Laplace, qu'il y en eut cinq \'editions entre 1814 et 1825. Cet essai philosophique est \`a la fois une merveilleuse introduction \`a la philosophie de Laplace, h\'erit\'ee des encyclop\'edistes, et au monumental trait\'e sur la th\'eorie analytique des probabilit\'es.

Pierre--Simon Laplace eut une longue vie : 1749--1827. C'est le cas pour d'autres grands math\'ematiciens de cette \'epoque : Euler et Lagrange qui le pr\'ec\`edent, Legendre qui fut exactement son contemporain, Gauss et Cauchy qui le suivent. Et ces math\'ematiciens furent actifs tout au cours de leur existence. Je le signale parce que la long\'evit\'e est fr\'equente chez les math\'ematiciens (Newton est mort \`a plus de 80 ans, Hadamard presque centenaire) et qu'une bonne part des \oe uvres est \'ecrite dans la maturit\'e et non dans la prime jeunesse. Les exemples fameux de math\'ematiciens morts jeunes (Galois, Abel) n'est sont que plus dramatiques.

\section*{Laplace acad\'emicien}

Laplace \'etait le fils d'un paysan pauvre. Il fut remarqu\'e, envoy\'e au Coll\`ege de Caen, puis \`a l'\'ecole militaire de Beaumont. Il y devint professeur \`a vingt ans. A vingt-deux ans il vint \`a Paris, rencontra d'Alembert, et obtint gr\^ace \`a lui une chaire \`a l'\'ecole militaire de Paris. L'histoire de la rencontre vaut d'\^etre connue. Laplace avait des recommandations pour d'Alembert, mais ne parvenait pas \`a le voir. Alors il s'avisa d'\'ecrire une lettre exposant ses id\'ees  sur les principes g\'en\'eraux de la m\'ecanique. Le jour m\^eme, d'Alembert le re\c cut.

Il fut nomm\'e membre  adjoint de l'Acad\'emie des Sciences \`a moins de vingt-cinq ans. L'Acad\'emie des Sciences de l'\'epoque, c'\'etait \`a la fois une soci\'et\'e savante et un petit CNRS. On pouvait y entrer jeune (Laplace n'avait \'ecrit que deux m\'emoires, spectaculaires il est vrai), et on y entrait pour y travailler. A partir de l\`a, il fut un personnage officiel et \'ecout\'e aussi bien sous la R\'evolution, sous l'Empire et sous la Restauration. Ces honneurs, autant que j'en puisse juger, n'\'ebranl\`erent pas son jugement. On cite de lui une r\'eponse fameuse \`a Napol\'eon, s'interrogeant sur le fait que Dieu et la cr\'eation n'apparaissaient nulle part dans son ouvrage de 1796, {\bf L'exposition du syst\`eme du monde} : ``Sire, je n'ai pas eu besoin de cette hypoth\`ese''.

C'est dans cet ouvrage --- \'ecrit sans faire usage d'aucun appareil math\'ema\-tique --- que Laplace \'emet une hypoth\`ese sur l'origine du syst\`eme solaire qui para\^{i}t encore aujourd'hui la plus vraisemblable :  n\'ebuleuse initiale, se condensant pour former le noyau central et les plan\`etes.

\section*{M\'ecanique c\'eleste et probabilit\'es}

Les deux sujets de pr\'edilection de Laplace furent la m\'ecanique c\'eleste et les probabilit\'es. Dans les deux domaines, il d\'eveloppa un merveilleux appareil math\'ematique, et l'appliqua \`a des situations concr\`etes. Les titres m\^emes \'evoquent ces deux aspects : {\bf Calcul int\'egral et syst\`eme du monde, Th\'eorie analytique des probabilit\'es.} Dans une vue assez dialectique de la r\'ealit\'e, qu'il exprima mieux que personne --- j'y reviendrai ---, il voyait \`a la fois le champ immense de la m\'ecanique issue de Newton (``l'\'etat pr\'esent de l'univers comme l'effet de son \'etat ant\'erieur, et comme la cause de celui qui va suivre''), et celui de la th\'eorie des hasards, qu'il pr\'ef\'era appeler th\'eorie des probabilit\'es (``presque toutes nos connaissances ne sont que probables ; et dans le petit nombre de choses que nous pouvons savoir avec certitude, dans les sciences math\'ematiques elles--m\^emes, les principaux moyens de parvenir \`a la v\'erit\'e, l'induction et l'analogie, se fondent sur les probabilit\'es ; en sorte que le syst\`eme entier des connaissances humaines se rattache \`a la th\'eorie des probabilit\'es''). Ainsi, tout au cours de son existence, il d\'eveloppa aussi bien les math\'ematiques du hasard et celles de la necessit\'e (probabilit\'es, et \'equations  diff\'erentielles). Il sut d\'evelopper les outils pour eux--m\^emes (fonctions g\'en\'eratrices, calcul symbolique) en les faisant fonctionner \`a l'int\'erieur des math\'ematiques pour retrouver et am\'eliorer les r\'esultats ant\'erieurs de l'analyse --- et il  ne manque pas une occasion de signaler ses sources, Wallis, Moivre, et bien d'autres. Il sut \'egalement les appliquer \`a une foule de probl\`emes et de situations : les {\bf  tables lunaires ou plan\'etaires} (o\`u, contrairement peut--\^etre \`a l'\'evidence, il applique non seulement la m\'ecanique c\'eleste, mais aussi les probabilit\'es), les {\bf anneaux de Saturne}, les {\bf mar\'ees}, les {\bf com\`etes, la rotation de la terre}, les {\bf g\'eod\'esiques} (probl\`eme pratique important au temps o\`u l'on mesurait le m\'eridien terrestre pour donner la d\'efinition du m\`etre), {\bf la pendule et la mesure du temps} la {\bf vitesse du son}, l'{\bf action de la lune sur l'atmosph\`ere}, la {\bf chaleur}, l'{\bf \'electricit\'e} (en collaboration avec Lavoisier), la {\bf d\'emographie}, et les {\bf statistiques}, les {\bf probabilit\'es des causes} et l'ensemble de ce qu'il appelait la {\bf philosophie naturelle}. J'ai cit\'e l\`a des titres ou des sujets de m\'emoires publi\'es entre 1774 et 1817. La plupart de ces travaux se trouvent embrass\'es dans deux trait\'es monumentaux : {\bf Le trait\'e de m\'ecanique c\'eleste}, publi\'e de 1799 \`a 1825, et {\bf La th\'eorie analytique des probabilit\'es}, publi\'es en 1812.

Laplace, hors ses d\'ebuts, ne semble pas avoir beaucoup enseign\'e. Mais ses \'ecrits font montre d'un souci  p\'edagogique constant. Trois \oe uvres sont \'ecrites dans une langue rigoureuse, sans aucun symbole math\'ematique (ainsi  $\pi$ se nomme le rapport de la circonf\'erence au diam\`etre), et elles m\'eriteraient d'\^etre des classiques de la litt\'erature fran\c caise. Il s'agit des {\bf dix le\c cons} aux Ecoles normales de 1795, de l'{\bf exposition du syst\`eme du monde} de 1796, et de l'{\bf essai philosophique sur les probabilit\'es}, de 1814--1825. Je dirai quelques mots des dix le\c cons, et je m'attarderai sur l'essai philosophique.

L'{\bf essai philosophique} d\'ebute par une phrase qui n'est pas insignifiante : ``Cet essai philosophique est le d\'eveloppement  d'une le\c con sur les probabilit\'es, que je donnai en 1795 aux Ecoles normales, o\`u je fus appel\'e comme professeur de math\'ematiques avec Lagrange par un d\'ecret de la Convention Nationale''. Ceci fut \'ecrit \`a la fin de l'Empire de Napol\'eon, \'edit\'e quatre fois au temps de la Restauration, sous Louis XVIII et Charles X, quand Laplace \'etait comte d'Empire puis marquis sous les deux rois. Laplace n'\'etait certes pas tenu d'\'evoquer l'\oe uvre de la Convention. Cette \'evocation dans la premi\`ere phrase de son essai me semble prouver deux choses : d'abord, que Laplace avait du caract\`ere ; ensuite, que la Convention avait compt\'e dans sa vie.

C'est une chose qu'on voit d\`es la premi\`ere des {\bf dix le\c cons}. Laplace y traite de la num\'eration. Il expose excellemment, la num\'eration d\'ecimale, la num\'eration binaire, la num\'eration duod\'ecimale, leurs avantages et inconv\'enients, puis il \'evoque bri\`evement mais de fa\c con  frappante les travaux de la commission constitu\'ee par la Convention pour adopter un syst\`eme unique, la raison du choix du syst\`eme d\'ecimal et du syst\`eme m\'etrique. Le\c con de math\'ematique --- parfaitement formalis\'ee, et en avance sur son temps --- m\^el\'ee \`a l'histoire du temps. Je n'examinerai pas plus avant le contenu des {\bf dix le\c cons} ; elles m\'eritent une \'etude \`a elles seules. L'{\bf exposition du syst\`eme du monde} en est un vaste compl\'ement, parce que Laplace n'avait pas eu le temps de toucher \`a la m\'ecanique c\'eleste, et, comme il l'\'ecrit, l'{\bf essai philosophique} est le prolongement de la dixi\`eme le\c con.

La conclusion de l'{\bf essai philosophique} se r\'ef\`ere au projet p\'edagogique de Laplace, et on verra par--l\`a qu'il \'etait en avance d'exactement cent cinquante ans sur les programmes de l'enseignement public en France, en ce qui concerne l'introduction des probabilit\'es :

\vskip4mm

\begin{minipage}{12cm}
\parindent=16pt
{\bf ``On voit par cet Essai que la th\'eorie des probabilit\'es  n'est au fond que le bon sens r\'eduit au calcul : elle fait appr\'ecier avec exactitude, ce que les esprits justes sentent d'instinct, sans qu'ils puissent souvent s'en rendre compte. Elle ne laisse rien d'arbitraire dans le choix des opinions et des partis \`a prendre, toutes les fois que l'on peut, \`a son moyen, d\'eterminer le choix le plus avantageux. Par l\`a, elle devient le suppl\'ement le plus heureux \`a l'ignorance et \`a la faiblesse de l'esprit humain. Si l'on consid\`ere les m\'ethodes  analytiques auxquelles cette th\'eorie a donn\'e naissance, la v\'erit\'e des principes qui lui servent de base, la logique fine et d\'elicate qu'exige leur emploi dans la solution des probl\`emes, les \'etablissements d'utilit\'e publique qui s'appuient sur elle, et l'extension qu'elle a re\c cu et qu'elle peut recevoir encore par son application aux questions les plus importantes de la Philosophie naturelle et des sciences morales ; si l'on observe ensuite que dans les choses m\^emes qui ne peuvent \^etre soumises au calcul elle donne les aper\c cus les plus s\^urs qui puissent nous guider dans nos jugements, et qu'elle apprend \`a se garantir des illusions qui souvent nous \'egarent, on verra qu'il n'est point de science plus digne de nos m\'editations, et qu'il soit plus utile de faire entrer dans le syst\`eme de l'instruction publique.''}
\end{minipage}

\section*{Laplace m\'ecaniste ?}

Revenons au d\'ebut de l'{\bf essai philosophique}. Ayant \'evoqu\'e les Ecoles normales et la Convention, puis son ouvrage sur la {th\'eorie analytique des probabilit\'es}, Laplace poursuit :

\vskip4mm
\begin{minipage}{12cm}
\parindent=16pt
{\bf ``Je pr\'esente ici, sans le secours de l'Analyse, les principes et les r\'esultats g\'en\'eraux de cette th\'eorie, en les appliquant aux questions les plus importantes de la vie, qui ne sont en effet, pour la plupart, que des probl\`emes de probabilit\'e. On peut m\^eme dire, \`a parler en rigueur, que presque toutes nos connaissances ne sont que probables ; et dans le petit nombre des choses que nous pouvons savoir avec certitude, dans les sciences math\'ematiques elles--m\^emes, les principaux moyens de parvenir \`a la v\'erit\'e, l'induction et l'analogie, se fondent sur les probabilit\'es ; en sorte que le syst\`eme entier des connaissances humaines se rattache \`a la th\'eorie expos\'ee dans cet essai.''}
\end{minipage}

\vskip4mm

En quelques pages alors, Laplace, dans une langue superbe, dit \`a la fois sa conception du d\'eterminisme, de la science, de la place qu'occupe la probabilit\'e, du d\'ebat d'opinion et de son r\^ole dans la soci\'et\'e. Ainsi, la fameuse profession de foi d\'eterministe fait partie d'un ensemble dont il ne faut pas la d\'etacher : la recherche de lois universelles, la confiance dans le d\'eveloppement de la science, la d\'enonciation de l'obscurantisme sous des formes vari\'ees, l'id\'ee tr\`es dialectique que la probabilit\'e est relative en partie \`a notre ignorance, en partie \`a nos connaissances, la relation entre probabilit\'e et certitude, l'analyse des erreurs et la tol\'erance pour les opinions diff\'erentes. Tout vaudrait d'\^etre cit\'e,  mais je me borne \`a quelques extraits. Il me para\^{i}t  important de souligner le titre  g\'en\'eral de ces quelques pages : {\bf De la probabilit\'e.}

Sur le d\'eterminisme et les lois universelles :

\vskip4mm

\begin{minipage}{12cm}
\parindent=16pt
{\bf ``Tous les \'ev\'enements, ceux m\^emes qui par leur petitesse semblent ne pas tenir aux grandes lois de la nature, en sont une suite aussi n\'ecessaire que les r\'evolutions du soleil.''...}
\end{minipage}
\eject
\begin{minipage}{12cm}
\parindent=16pt
{\bf  ``Nous devons donc envisager l'\'etat pr\'esent de l'univers comme l'effet de son \'etat ant\'erieur, et comme la cause de celui qui suivre. Une intelligence qui pour un instant donn\'e conna\^{i}trait toutes les forces dont la nature est  anim\'ee et la situation respective des \^etres qui la composent, si d'ailleurs elle \'etait assez vaste pour soumettre ces donn\'ees \`a l'analyse, embrasserait dans la m\^eme formule les mouvements des plus grands corps de l'univers et ceux du plus l\'eger atome : rien ne serait incertain pour elle, et l'avenir comme le pass\'e serait pr\'esent \`a ses yeux.''}
\end{minipage}

\vskip4mm

\dots\  la confiance dans l'esprit humain :

\vskip4mm

\begin{minipage}{12cm}
\parindent=16pt
{\bf ``L'esprit humain offre, dans la perfection qu'il a su donner \`a l'Astronomie, une faible esquisse de cette  intelligence. Ses d\'ecouvertes en m\'ecanique et en g\'eom\'etrie, jointes \`a celle de la pesanteur universelle, l'ont mis \`a port\'ee de comprendre dans les m\^emes expressions analytiques les \'etats pass\'es et futurs du syst\`eme du monde. En appliquant la m\^eme m\'ethode \`a quelques autres objets de ses connaissances , il est parvenu \`a ramener \`a des lois g\'en\'erales les ph\'enom\`enes observ\'es, et \`a pr\'evoir ceux que des circonstances donn\'ees doivent faire \'eclore. Tous ces efforts dans la recherche de la  v\'erit\'e tendent \`a le rapprocher sans cesse de l'intelligence que nous venons  de concevoir, mais dont il restera toujours infiniment \'eloign\'e. Cette tendance propre \`a l'esp\`ece humaine est ce qui la rend sup\'erieure aux animaux, et  ses progr\`es en ce genre distinguent les nations et les si\`ecles et font leur v\'eritable gloire.''}
\end{minipage}

\vskip4mm

\dots \  la d\'enonciation des superstitions :

\vskip4mm

\begin{minipage}{12cm}
\parindent=16pt
{\bf ``Rappelons-nous qu'autrefois, et  \`a une  \'epoque qui n'est pas encore bien recul\'ee, une pluie ou une s\'echeresse extr\^eme,  une com\`ete tra\^{i}nant apr\`es elle une queue fort \'etendue, les \'eclipses, les aurores bor\'eales, et g\'en\'eralement tous les ph\'enom\`enes extraordinaires, \'etaient regard\'es\break}
\end{minipage}
\eject
\begin{minipage}{12cm}
{\bf comme autant de signes de la col\`ere c\'eleste. On invoquait le ciel pour d\'etourner leur funeste influence. On ne le priait point de suspendre le cours des plan\`etes et du soleil : l'observation e\^ut bient\^ot fait sentir  l'inutilit\'e de ces pri\`eres. Mais comme ces ph\'enom\`enes, arrivant et disparaissant \`a de longs intervalles, semblaient contrarier l'ordre de la nature, on supposait que le ciel, irrit\'e par les crimes de la terre, les faisait  na\^{i}tre pour annoncer ses vengeances. Ainsi la longue queue de la com\`etre de 1456 r\'epandit la terreur dans l'Europe, d\'ej\`a constern\'ee par les succ\`es rapides des Turcs qui venaient de renverser le Bas--Empire. Cet astre, apr\`es quatre de ses r\'evolutions, a excit\'e parmi nous un int\'er\^et bien diff\'erent. La connaissance des lois du syst\`eme du monde, acquise dans cet intervalle, avait dissip\'e les craintes enfant\'ees par l'ignorance des vrais rapports de l'homme avec l'univers ; et Halley, ayant reconnu l'identit\'e de cette com\`ete avec celles des ann\'ees 1531, 1607 et 1682, annon\c ca son retour prochain pour la fin de 1758 ou le commencement de 1759. Le monde savant attendit avec impatience ce retour qui devait confirmer l'une des plus grandes d\'ecouvertes que l'on e\^ut faites dans les sciences,''\dots}
\end{minipage}

\vskip4mm

\dots\ L'analyse et les lois universelles :

\vskip4mm

\begin{minipage}{12cm}
\parindent=16pt
{\bf ``Clairaut entreprit alors de soumettre \`a l'analyse les perturbations que la com\`ete avait par l'action des deux plus grosses plan\`etes, Jupiter et Saturne : apr\`es d'immenses calculs, il fixa son prochain passage au p\'erih\'elie vers le commencement d'avril 1759, ce que l'observation ne tarda pas \`a v\'erifier. La r\'egularit\'e que l'Astronomie nous montre dans le mouvement des com\`etes a lieu, sans aucun doute, dans tous les ph\'enom\`enes. La courbe d\'ecrite par une simple mol\'ecule d'air ou de vapeurs est r\'egl\'ee d'une mani\`ere  aussi certaine que les orbites plan\'etaires : il n'y a pas de diff\'erence entre elles que celle qu'y met notre ignorance''}
\end{minipage}

\eject

\dots\ La probabilit\'e :

\vskip4mm

\begin{minipage}{12cm}
\parindent=16pt

{\bf ``La probabilit\'e est relative en partie \`a cette ignorance, en partie \`a nos connaissances\dots''

``La th\'eorie des hasards consiste \`a r\'eduire tous les \'ev\'enements du m\^eme genre \`a un certain nombre de cas \'egalement possibles, c'est--\`a--dire tels que nous soyons \'egalement ind\'ecis sur leur existence, et \`a d\'eterminer le nombre de cas favorables \`a l'\'ev\'enement dont on cherche la probabilit\'e. Le rapport de ce nombre \`a tous les cas possibles est la mesure de cette probabilit\'e qui n'est ainsi qu'une fraction, dont le num\'erateur est le nombre des cas favorables et dont le d\'enominateur est le nombre de tous les cas possibles.''\dots}
\end{minipage}

\vskip4mm

(Je me permets d'insister ici sur la pr\'ecision du langage de Laplace, qui parle de cas {\bf \'egalement possibles} ; nous dirions aujourd'hui, moins \'el\'egamment, \'equiprobables).

\vskip4mm
\begin{minipage}{12cm}
\parindent=16pt

{\bf ``Quand tous les cas sont favorables \`a un \'ev\'enement, sa probabilit\'e se change en certitude et son expression devient \'egale \`a l'unit\'e. Sous ce rapport, la certitude et la probabilit\'e sont comparables, quoiqu'il y ait une diff\'erence essentielle, entre les deux \'etats d'esprit, lorsqu'une v\'erit\'e lui est rigoureusement d\'emontr\'ee, ou lorsqu'il aper\c coit encore une petite source d'erreur.''}
\end{minipage}

\vskip4mm

\dots\ la diff\'erence des opinions, les erreurs, le pluralisme :

\vskip4mm

\begin{minipage}{12cm}
\parindent=16pt

{\bf ``Dans les choses qui ne sont que vraisemblables, la diff\'erence des donn\'ees que chaque homme a sur elles est une des causes principales de la diversit\'e des opinions que l'on voit r\'egner sur les m\^emes objets.''\dots

``C'est \`a l'influence de l'opinion de ceux que la multitude juge les plus instruits, et \`a qui elle a coutume de donner sa confiance sur les plus importants objets de la vie, qu'est due la propagation de ces erreurs qui, dans les\break}
\end{minipage}
\eject

\begin{minipage}{12cm}
{\bf temps d'ignorance, ont couvert la face du monde. La Magie et l'Astrologie nous en offrent deux grands exemples. Ces erreurs inculqu\'ees d\`es l'enfance, adopt\'ees sans examen, et n'ayant pour base que la croyance universelle, se sont maintenues pendant tr\`es longtemps ; jusqu'\`a ce qu'enfin le progr\`es des sciences les ait  d\'etruites dans l'esprit des hommes \'eclair\'es dont ensuite l'opinion les a fait dispara\^{i}tre chez le peuple m\^eme, par le pouvoir de l'imitation et de l'habitude, qui les avait si g\'en\'eralement r\'epandues. Ce pouvoir, le plus puissant ressort du monde moral, \'etablit et conserve dans toute une nation des id\'ees enti\`erement contraires \`a celles qu'il maintient ailleurs avec le m\^eme empire. Quelle indulgence ne devons--nous donc pas avoir pour les opinions  diff\'erentes des n\^otres, puisque cette diff\'erence ne d\'epend souvent que des points divers o\`u les circonstances nous ont plac\'es ! Eclairons ceux que nous ne jugeons pas suffisamment instruits, mais auparavant examinons s\'ev\`erement nos propres opinions, et pesons avec impartialit\'e leurs probabilit\'es respectives.''}
\end{minipage}

\vskip4mm

Apr\`es ces quelques citations, on conviendra qu'il est un peu sch\'ematique de repr\'esenter Laplace comme un m\'ecaniste \`a tout crin, soucieux seulement de d\'emonter l'univers comme une grande horloge. Au contraire, une grande partie de son \oe uvre, et l'essentiel de l'{\bf essai philosophique}, concerne les relations entre le contingent, le possible, le probable  et le certain. Devant le hasard nous ne somme pas d\'emunis : en formalisant le peu que nous connaissons, en construisant des mod\`eles convenables, en appliquant l'analyse math\'ematique, nous pouvons arriver \`a ces conclusions nettes en face de ph\'enom\`enes al\'eatoires.

\section*{Laplace th\'eoricien des probabilit\'es}

L'essai compte plus de deux cents pages. Je me bornerai maintenant \`a un survol, avec quelques plong\'ees pour donner un certain go\^ut du contenu.

\subsection*{Le passage du fini \`a l'infini}

Le premier quart, relatif aux principes  et aux m\'ethodes, est un tour de force. En dix principes, il explique le fondement de la th\'eorie des probabilit\'es, par des d\'efinitions impeccables, assorties d'exemples. La terminologie n'a gu\`ere vari\'e depuis Laplace (sinon que ``croix ou pile'' s'appelle maintenant ``pile ou face''). La th\'eorie moderne de la mesure, d\'evelopp\'ee depuis Lebesgue (1904), permet d'aller plus vite et plus loin, et les expos\'es de base que nous faisons aujourd'hui s'inspirent plus ou moins de l'axiomatique de Kolmogorov (1933). Laplace d\'eveloppe ses principes sur des mod\`eles finis (urnes, jeux de pile ou face), tandis que l'essentiel de son  \oe uvre math\'ematique concerne le passage du fini \`a l'infini.  Il n'y a pas l\`a de vice logique. Pour lui, les lois limites  (telle la fameuse loi de Gauss, dont il fait la th\'eorie et l'utilisation quand Gauss \'etait petit enfant !), sont des moyens analytiques pour approcher les probabilit\'es relatives \`a un mod\`ele fini.

Ces m\'ethodes analytiques sont expos\'ees \`a la suite. La m\'ethode des fonctions  g\'en\'eratrices, issue de probl\`emes de probabilit\'e, s'applique \'egalement dans d'autres domaines (\'equations aux diff\'erences, \'equations diff\'erentielles). Laplace parvient \`a l'exposer sans \'ecrire une formule. Voici maintenant comment il introduit la loi de Laplace--Gauss :

\vskip4mm
\begin{minipage}{12cm}
\parindent=16pt

{\bf ``On est souvent conduit \`a des expressions qui contiennent tant de termes et de facteurs, que les substitutions num\'eriques y sont impraticables. C'est ce qui a lieu dans les questions de probabilit\'e, lorsque l'on consid\`ere un grand nombre d'\'ev\'enements. Cependant, il importe alors d'avoir la valeur num\'erique des formules, pour conna\^{i}tre avec quelle probabilit\'e les r\'esultats, que les \'ev\'enements d\'eveloppent en se multipliant, sont indiqu\'es. Il  importe surtout d'avoir la loi suivant laquelle cette probabilit\'e approche sans cesse de la certitude  qu'elle finirait par atteindre, si le nombre des \'ev\'enements devenait infini.
Pour y parvenir, je consid\'ererai que les int\'egrales d\'efinies de diff\'erentielles multipli\'ees par des facteurs \'elev\'es \`a de grandes puissances, donnaient par l'int\'egration des formules\break}
\end{minipage}
\eject

\begin{minipage}{12cm}
{\bf  compos\'ees d'un grand nombre de termes et de facteurs. Cette remarque me fit na\^{i}tre  l'id\'ee de transformer dans de semblables int\'egrales les expresssions compliqu\'ees de l'analyse, et les int\'egrales des \'equations aux diff\'erences. J'ai rempli cet objet, par une m\'ethode qui donne \`a la fois la fonction comprise sous le signe int\'egral, et les limites de l'int\'egration.''}
\end{minipage}

\vskip4mm

Puis il indique le lien avec la th\'eorie des fonctions g\'en\'eratrices, et avec les travaux de ses pr\'ed\'ecesseurs (avec des aper\c cus pr\'ecieux sur le r\^ole des notations dans la pens\'ee math\'ematique, l'influence de Fermat, de Descartes, de Wallis, de Newton, de Leibniz, de Jacques Bernoulli, de Lagrange, etc\dots)

Presque tout le reste de l'{\bf essai philosophique} est consacr\'e aux applications. Cela touche \`a tout : aux naissances de gar\c cons et de filles, au calcul des erreurs, \`a la correction des donn\'ees astronomiques, aux op\'erations g\'eod\'esiques (mesure du m\'eridien terrestre), aux in\'egalit\'es dans les mouvements de la lune, \`a l'applatissement de la Terre, aux mouvements de Jupiter et de Saturne et \`a ceux des satellites de Jupiter, aux mar\'ees, \`a l'influence de la lune sur l'atmosph\`ere, aux variations de presssion atmosph\'erique et de magn\'etisme, \`a la th\'eorie du syst\`eme solaire, aux com\`etes, aux t\'emoignages en justice, aux choix et d\'ecisions des assembl\'ees, aux assurances, \`a la d\'emographie, aux maladies et  \`a la vaccination, etc.

Enfin, une quarantaine de pages traitent des illusions dans l'estimation des probabilit\'es : du r\^ole des habitudes, des croyances, de la m\'emoire, des m\oe urs, des exp\'eriences en psychologie, des somnambules, des visionnaires, de la musique, de l'attention, des images mentales, etc. Il ne serait pas sans int\'er\^et de lire ces pages \`a la lumi\`ere de la neurophysiologie moderne. Elles  contiennent des vues justes et int\'eressantes. Cependant, que d'\'ecart entre le ``sensorium''  de Laplace, objet de vibrations soumises aux lois de la dynamique, et l'``homme neuronal''  de Changeux !

\subsection*{La statistique}

A part cette incursion, parfois hasardeuse, dans le domaine de l'imaginaire, rien n'a vieilli. Je voudrais le montrer  sur quelques exemples.

Laplace  --- le fait est insuffisamment connu --- est non seulement le premier th\'eoricien des probabilit\'es, mais le fondateur de la statistique math\'ematique. Les ``fourchettes''  et les tests d'hypoth\`ese ont \'et\'e introduits \`a l'occasion d'une question bizarre, celle de la proportion de gar\c cons et de filles \`a la naissance. Cette proportion est  ordinairement de 22/21 (noter qu'on utilisait encore rarement les nombres d\'ecimaux \`a cette \'epoque, au moins dans la langue courante). On a signal\'e qu'elle \'etait invers\'ee, dans un bourg de Bourgogne, Carcelle--les--Grignon, pendant 5 ans : sur 2~009 naissances, 1~026 filles et 983 gar\c cons. Sur ces chiffres, Laplace teste l'hypoth\`ese admise, \`a savoir que les naissances de gar\c cons sont plus  probables que celles des filles ; le calcul (n\'ecessitant l'approximation par la fonction de Laplace--Gauss) conclut \`a ne pas rejeter l'hypoth\`ese. En revanche, \`a Paris, sur une p\'eriode de quarante ans, on recense \`a la naissance 393~386 gar\c cons et 377~555 filles ; la proportion est de 25/24, assez proche de 22/21. Laplace teste alors l'hypoth\`ese  qu'\`a Paris les probabilit\'es de naissances masculines et f\'eminines soient proportionnelles \`a 22 et 21. L\`a  il rejette cette hypoth\`ese. Il faut donc une explication du fait qu'\`a Paris on recense un peu moins de gar\c cons qu'en g\'en\'eral. Laplace trouve un \'el\'ement  de r\'eponse. A Paris plus de 20~\% des enfants sont des enfants trouv\'es. ``En y r\'efl\'echissant, dit Laplace, il m'a paru que la diff\'erence observ\'ee tient \`a ce que les parents de la campagne et des provinces, trouvant quelque avantage \`a retenir pr\`es d'eux  les gar\c cons, en avaient envoy\'e \`a l'Hospice des enfants trouv\'es de Paris moins relativement aux filles, que suivant le rapport des naissances des deux sexes''. V\'erification : Laplace soustrait des  effectifs parisiens ceux des enfants trouv\'es, et trouve exactement 22/21.

Dans la derni\`ere \'edition de son essai, Laplace est troubl\'e. Il donne les tableaux des naissances entre 1817 et 1821 (les cinq premi\`eres ann\'ees d'un recensement g\'en\'eral des naissances en France) : 2~962~361 gar\c cons, 2~781~997 filles, proportion 16/15. L'\'ecart 22/21, compte tenu des nombres en jeu, rend inadmissible l'hypoth\`ese d'un rapport constant, et Laplace n'a pas d'explication. Nous n'avons pas non plus, mais le fait est maintenu reconnu que ce rapport change apr\`es les grandes guerres (en France, par exemple, apr\`es 1918). Laplace, en effet, ne trichait pas avec le calcul ;  j'en indiquerai tout \`a l'heure un autre exemple.

\subsection*{Les enqu\^etes}

La statistique s'est impos\'ee \`a Laplace pour des raisons moins gratuites. Il fut ministre de l'Int\'erieur apr\`es le 18~ brumaire, et il fallait conna\^{i}tre la population, d\'epartement par d\'epartement, pour la conscription. En l'absence de recensement g\'en\'eral, il pratiqua un  v\'eritable plan d'exp\'erience, en utilisant la proportionnalit\'e suppos\'ee entre la population et le nombre des naissances. Ce fut la premi\`ere estimation donn\'ee de la population fran\c caise : 28~352~845, \`a un demi--million pr\`es, avec un risque d'erreur inf\'erieur \`a 1/300~000. Laplace d\'ecrit ainsi le plan pour d\'eterminer le  rapport de proportionnalit\'e entre la population et le nombre des naissances : il consiste

\vskip4mm

\begin{minipage}{12cm}
\parindent=16pt

{\bf ``1) \`a choisir dans l'empire les d\'epartements distribu\'es d'une mani\`ere \`a peu pr\`es \'egale sur toute sa surface, afin de rendre le r\'esultat g\'en\'eral ind\'ependant des circonstances locales ;
2) \`a d\'enombrer avec soin, pour une  \'epoque donn\'ee, les habitants de plusieurs communes dans chacun de ces d\'epartements ; 3) \`a d\'eterminer par le relev\'e des naissances durant plusieurs ann\'ees, qui pr\'ec\`edent et suivent cette \'epoque, le nombre moyen correspondant des naissances annuelles.''\dots}
\end{minipage}

\vskip4mm

Et il indique comment, sur sa pri\`ere, le gouvernement a fait mener l'enque\^ete, dans trente d\'epartements ``r\'epandus \'egalement dans toute la France'', en faisant choix des communes pouvant donner les renseignements les plus pr\'ecis : en 1800, 1801, 1802, 116~599 naissances ; au 23 septembre, 2~037~615 habitants ; rapport :  28,352~845 (cela correspondait \`a l'esp\'erance de vie d'un nouveau n\'e).

L'autorit\'e de Laplace et ses responsabilit\'es officielles l'ont amen\'e \`a demander et obtenir d'autres statistiques. C'est ainsi qu'\`a dater du 1\up{er} juin 1806 le niveau des mar\'ees a \'et\'e not\'e, chaque jour, dans  le port de Brest. Il  s'agissait de tester la th\'eorie de Laplace sur les mar\'ees, un immense travail qui conclut \`a la validit\'e parfaite des formules fond\'ees sur l'influence de la lune et du soleil.

Laplace a \'egalement propos\'e des formules pour des mouvements \`de l'atmosph\`ere sous l'action de la lune et du soleil, et il a pu prescrire huit ann\'ees d'observations de la pression atmosph\'erique \`a Paris quatre fois par jour. Il observe une faible oscillation, \`a l'appui de son hypoth\`ese. La suite m\'erite d'\^etre cit\'ee comme un mod\`ele d'honn\^etet\'e scientifique.

\vskip4mm

\begin{minipage}{12cm}
\parindent=16pt

{\bf ``C'est ici surtout que se fait sentir la n\'ecessit\'e d'une m\'ethode pour d\'eterminer la probabilit\'e d'un r\'esultat, m\'ethode sans laquelle on est expos\'e \`a pr\'esenter comme lois de la nature, les effets des causes irr\'eguli\`eres, ce qui est arriv\'e souvent en M\'et\'eorologie.\break}
\end{minipage}
\eject
\begin{minipage}{12cm}
{\bf
 Cette m\'ethode appliqu\'ee au r\'esultat pr\'ec\'edent en montre l'incertitude, malgr\'e le grand nombre d'observations employ\'ees, qu'il faudrait d\'ecupler pour obtenir un r\'esultat suffisamment probable.''}
\end{minipage}

\vskip4mm

En voici un autre exemple, sur un sujet qu'il n'a pas \'etudi\'e, mais o\`u il d\'efinit clairement la bonne attitude \`a prendre

\vskip4mm

\begin{minipage}{12cm}
\parindent=16pt

{\bf ``Les ph\'enom\`enes singuliers qui r\'esultent de l'extr\^eme sensibilit\'e des nerfs dans quelques individus, ont donn\'e naissance \`a diverses opinions sur l'existence d'un nouvel agent que l'on a nomm\'e magn\'etisme animal, sur l'action du magn\'etisme ordinaire, et sur l'influence du soleil et de la lune dans quelques affections nerveuses, enfin sur les impresssions que peut faire \'eprouver la proximit\'e des m\'etaux ou d'une eau courante. Il est naturel de penser que l'action de ces causes est tr\`es faible, et qu'elle peut \^etre facilement troubl\'ee par des circonstances accidentelles ; ainsi, parce que dans quelques cas elle ne s'est point manifest\'ee, on ne doit pas rejeter son existence. Nous sommes si loin de conna\^{i}tre tous les agents de la nature et leurs divers modes d'action, qu'il serait peu philosophique de nier les ph\'enom\`enes, uniquement parce qu'ils sont inexplicables dans l'\'etat actuel de nos connaissances. Seulement, nous devons les examiner avec une attention d'autant plus scrupuleuse, qu'il para\^{i}t  plus difficile de les admettre ; et c'est ainsi que le calcul des probabilit\'es devient indispensable, pour d\'eterminer jusqu'\`a quel point il faut multiplier les observations ou les exp\'eriences, afin d'obtenir en faveur des agents qu'elles indiquent une probabilit\'e sup\'erieure aux raisons que l'on peut  avoir d'ailleurs, de ne pas les admettre.''}
\end{minipage}

\subsection*{Les sciences morales}

Laplace consacre cinquante pages de son {\bf essai} aux applications des probabilit\'es aux sciences morales. Faute de les analyser, je veux en indiquer le grand int\'er\^et (\`a la fois pour le sujet lui--m\^eme,  et par ce qui se r\'ev\`ele de la personnalit\'e de Laplace). Au plan scientifique, il est remarquable qu'il pose clairement, d\`es  cette \'epoque, la question de la mod\'elisation math\'ematique, \`a propos d'une question difficile, la probablit\'e des t\'emoignages.

\vskip4mm

\begin{minipage}{12cm}
\parindent=16pt

{\bf ``La plupart de nos jugements \'etant fond\'es sur la probabilit\'e des t\'emoignages, il est bien important de la soumettre au calcul. La chose, il est vrai, devient souvent impossible, par la difficult\'e d'appr\'ecier la v\'eracit\'e des t\'emoins et par le grand nombre de circonstances dont les faits qu'ils attestent sont accompagn\'es. Mais on peut dans plusieurs cas r\'esoudre  des probl\`emes qui ont beaucoup d'analogies avec les questions qu'on se propose, et  dont les solutions peuvent \^etre regard\'ees comme des approximations propres \`a nous guider et \`a nous  garantir des erreurs et des dangers auxquels de mauvais raisonnements nous exposent.''}
\end{minipage}

\section*{Laplace \'ecrivain}

Nous avons la chance en France d'avoir eu de grands savants qui ont \'et\'e de grands \'ecrivains. Laplace est de ceux--l\`a, j'esp\`ere vous en avoir convaincus. Tirons--nous suffisamment parti ce cette richesse ? Cherchons--nous suffisamment, chacun \`a notre niveau, de la prendre en exemple dans notre propre travail ? Les r\'eponses \`a ces deux questions sont \'evidemment n\'egatives, et il serait dommage de donner \`a cet article une conclusion formellement n\'egative. Pr\'esentons les questions autrement : pouvons--nous et devons--nous faire mieux ? Oui, assur\'ement.

\vskip2mm

\hfill Jean--Pierre Kahane

\hfill 01.08.07

Cet article recoupe largement un appel \`a r\'e\'editer l'Essai philosophique sur les probabilit\'es, paru dans l'Ecole et la Nation n\up{o} 351 (ao\^ut--septembre 1984), pp.~40--44.

\end{document}